\documentclass[12pt]{amsart}

\usepackage{amsmath, amssymb, latexsym, verbatim}

\newtheorem{Thm}{Theorem}[section]
\newtheorem{Lem}{Lemma}[section]

\newtheorem{Cor}{Corollary}[section]
\newtheorem{Def}{Definition}[section]

\def\H{\mathbb{H}}
\def\R{\mathbb{R}}
\def\SS{\mathbb{S}}
\def\N{\mathbb{N}}
\def\L{\mathbb{L}}

\def\bX{\bar{X}}
\def\bg{\bar{g}}
\def\gH{g_{_{\mathbb{H}}}}
\def\hg{\hat{g}}
\def\l({\left(}
\def\r){\right)}

\begin{document}

\title{Hypersurfaces in Hyperbolic Poincar\'{e} Manifolds and
Conformally Invariant PDEs}

\date{October 2009}


\renewcommand{\subjclassname}{\textup{2000} Mathematics Subject Classification}

\author{Vincent Bonini \, Jos\'{e} M. Espinar and Jie Qing }

\address{Department of Mathematics, California Polytechnic State
University, San Luis Obispo, CA 93407.} \email{vbonini@calpoly.edu}

\address{Departmento de Geometr\'{i}a y Topolog\'{i}a, Universidad de
Granada, E-18071 Granda, Spain.}
\email{jespinar@ugr.es}

\address{Department of Mathematics, University of California, Santa
Cruz, CA 95064.} \email{qing@ucsc.edu}

\begin{abstract}
We derive a relationship between the eigenvalues of the
Weyl-Schouten tensor of a  conformal representative of the conformal
infinity of a hyperbolic Poincar\'{e} manifold and the principal
curvatures on the level sets of its uniquely associated defining
function with calculations based on \cite{FG1} \cite{FG2}. This
relationship generalizes the result for hypersurfaces in
${\H}^{n+1}$ and their connection to the conformal geometry of
${\SS}^n$ as exhibited in \cite{JE} and gives a correspondence
between Weingarten hypersurfaces in hyperbolic Poincar\'{e}
manifolds and conformally invariant equations on the conformal
infinity. In particular, we generalize an equivalence exhibited in
\cite{JE} between Christoffel-type problems for hypersurfaces in
${\H}^{n+1}$ and scalar curvature problems on the conformal infinity
${\SS}^n$ to hyperbolic Poincar\'{e} manifolds.
\end{abstract}

\maketitle

\markboth{Bonini, Espinar and Qing}{}

\section{Introduction}

The relationship between the geometry of a conformally compact space
and the geometry of its conformal infinity has been of recent
interest in both physical and mathematical communities. The interest
in such association is motivated primarily by the AdS/CFT
correspondence where a conformal field theory on a compact manifold
$M^{n}$ correlates to the string theory on a negatively curved
conformally compact Einstein manifold $X^{n+1}$, which has $M$ as
its conformal infinity. One can view such connections as originating from the
identification between the group of isometries of hyperbolic space
${\H}^{n+1}$ and the group of conformal transformations of the round
sphere ${\SS}^n$. In fact, the study of such connections date
back to the 1980's in the seminal paper of Fefferman and Graham
\cite{FG1}.

Recently, an explicit example connecting the geometry of hyperbolic
space ${\H}^{n+1}$ to the conformal geometry of the round sphere
${\SS}^n$ was realized in \cite{JE} in the context of the
hypersurface geometry of ${\H}^{n+1}$ and curvature prescription
problems on ${\SS}^n$ in the conformal class of the round metric.
The question of the existence a hypersurface $\Sigma^n$ in
hyperbolic space ${\H}^{n+1}$ with prescribed Weingarten functional
of the principal curvatures of $\Sigma$ is a natural extension of
the classical problem in the Euclidean setting. Of particular
interest is the Christoffel problem for hypersurfaces in hyperbolic
space ${\H}^{n+1}$ where one is asked to find a hypersurface
$\Sigma^{n} \subset {\H}^{n+1}$ with prescribed mean of the
curvature radii. One of the initial difficulties of the Christoffel
problem in ${\H}^{n+1}$ is to provide the appropriate formulation of
the Gauss map and the principal curvature radii in the context of
hyperbolic space. In \cite{JE} the relevant notions of the
hyperbolic Gauss map and the hyperbolic principal curvature radii
are developed using the ambient structure of the hyperboloid model
of ${\H}^{n+1}$ where hyperbolic space is realized as a hypersurface
in Minkowski spacetime. Moreover, \cite{JE} exhibits a strikingly
precise relation between Christoffel-type problems for immersed
hypersurfaces in ${\H}^{n+1}$ and scalar curvature prescription
problems of conformal geometry on ${\SS}^n$, viewed as the boundary
of ${\H}^{n+1}$ at infinity. See also a related work of Mazzeo and
Pacard \cite{MP}.

In this note we take a viewpoint more reflective of conformal
geometry and we generalize the correspondences exhibited in
\cite{JE} between Christoffel-type problems and scalar curvature prescription
problems of conformal geometry. For $n \geq 2$, let $X^{n+1}$ denote the
interior of a smooth compact manifold ${\bX}^{n+1}$ with boundary
$\partial X = M^n$. A Riemannian metric $g$ on $X$ is then said to be
conformally compact if, for a defining function $r$ for $M$, the conformal
metric $\bg = r^2 g$ extends to a metric on $\bX$. The metric $\bg$ restricted
to $T M$ induces a metric $\hg$ on $M$, which rescales by conformal factor
upon change in defining function and therefore defines a conformal structure
$(M, [\hg])$ on $M$ called the conformal infinity of $(X, g)$.

A {\em hyperbolic Poincar\'{e} manifold} is a conformally compact
hyperbolic manifold. From the work of \cite{FG2}, given a representative
$\gamma\in [\hg]$ of the conformal infinity of a hyperbolic Poincar\'{e} manifold
$(X^{n+1}, g)$ and its associated geodesic defining function $r$, we may
write the metric in the normal form
\begin{equation} \label{E: FG Norm Form}
g = r^{-2}( dr^2 + g_r)
\end{equation}
where
\begin{equation}\label{E: FG}
g_r = \gamma - r^2P_{\gamma} + \frac{r^4}{4} Q(P_{\gamma}),
\end{equation}
$$
Q(P_\gamma)_{ij} = \gamma^{kl}(P_\gamma)_{ik}(P_\gamma)_{jl}
$$
and for $n \geq 3$,
$$
P_\gamma = \frac 1{n-2}\l( Ric_\gamma - \frac
{R_\gamma}{2(n-1)}\gamma\r)
$$
is the Schouten tensor of $\gamma$ with $Ric_\gamma$ and $R_\gamma$
denoting the Ricci and the scalar curvature of $\gamma$, respectively (please refer to
\S\ref{S: AH Mnflds} for more details). For $n=2$, $P_\gamma$ is a symmetric
2-tensor on $M$ satisfying
$$
\gamma^{ij} (P_\gamma)_{ij}  = \frac{ R_{\gamma}}{2} \quad
\text{and} \quad \gamma^{jk} (P_\gamma)_{ij,k}= (R_{\gamma})_{,i}.
$$

We will show that the horospherical metric associated to a
horospherical ovaloid in hyperbolic space $\H^{n+1}$ can be realized
as a representative of the conformal infinity $(\SS^n, [g_0])$ of
hyperbolic space $(\H^{n+1}, \gH)$. This is because a horospherical
ovaloid in hyperbolic space $\H^{n+1}$ determines a geodesic
defining function $r$ for the infinity $\SS^n$ of $\H^{n+1}$, where
$r = e^{- s}$ and $s$ is the hyperbolic distance to the
horospherical ovaloid. In general, on an asymptotically hyperbolic
manifold $X$, one should replace the notion of a horospherical
ovaloid by that of an {\em essential set}. As defined in \cite{BM},
the exponential map from the normal bundle of an essential set is a
diffeomorphism to the outside of the essential set in $X$. Hence, an
essential set provides a geodesic defining function $r = e^{-s}$
where $s$ is the distance to the essential set in $X$. A similar
idea was realized in early works of Epstein \cite{Epstein1}
\cite{Epstein2}.

A straightforward calculation based on (\ref{E: FG}) yields a
generalization of the relation in \cite{JE} between the eigenvalues
of the Schouten tensor of the horospherical metric and the
hyperbolic principal curvature radii of the level sets of the
associated geodesic defining function. To avoid any possible sign
confusion of the principal curvature of a hypersurface we recall
that, the second fundamental of a hypersurface $\Sigma$ in $X^{n+1}$
with respect to an orientation induced by a choice a normal
direction $N$ is defined to be
\begin{equation}\label{secondfundamentalform}
II = - \frac 12 \mathcal {L}_N g,
\end{equation}
where $\mathcal {L}$ is the Lie derivative. In our convention, for
instance, the principal curvature of a unit sphere in Euclidean
space with the orientation induced by the inward normal direction is
$1$.

\begin{Thm} \label{T: Main Thm}
Suppose that $(X^{n+1}, g)$ is a hyperbolic Poincar\'{e} manifold
and let $\gamma$ be a representative of its conformal infinity $(M^n, [\hg])$
with associated geodesic defining function $r$. Then the eigenvalues
$\lambda_i$ of the tensor $P_{\gamma}$ in the expansion (\ref{E: FG}) satisfy
\begin{equation}\label{E: key}
1 - \frac{r^2}{2} \lambda_i = \frac{2}{1-\kappa_i}
\end{equation}
where $\kappa_i=\kappa_i(r)$ denotes the $i^{th}$ outward principal
curvature on the level sets of the geodesic defining function $r$
and $\frac{2}{1- \kappa_i}$ is considered to be the $i^{th}$
hyperbolic principal curvature radius.
\end{Thm}

As studied in \cite{JE}, when $n\geq 3$, the relationship (\ref{E: key}) in Theorem \ref{T: Main Thm}
can be used to turn questions regarding foliations near the conformal infinity by particular classes of
hypersurfaces in hyperbolic Poincar\'{e} manifolds into questions regarding the conformal
geometry of the conformal infinity and visa versa. For example, taking the trace of (\ref{E: key}),
one finds that
\begin{equation} \label{E: Scal Curv Corr}
R_{\gamma} = \frac{4(n-1)}{r^2} \l(n - \sum_{i=1}^{n} \frac{2}{1 -
\kappa_i}\r).
\end{equation}
Therefore, finding a foliation by hypersurfaces with constant mean
of the hyperbolic curvature radii is equivalent to finding a
constant scalar curvature metric on the conformal infinity. Hence,
due to the resolution of the Yamabe problem we have the following Corollary.

\begin{Cor}\label{C: main}
Suppose that $(X^{n+1}, g)$ is a hyperbolic Poincar\'{e} manifold.
Then there always exists a foliation of hypersurfaces of constant
mean of the hyperbolic curvature radii near the infinity.
Such foliations are parameterized by geodesic defining functions $r$
associated with constant scalar curvature $S$ representatives of the
conformal infinity and the mean of the hyperbolic curvature radii of the
foliation is given by
\begin{equation} \label{E: Mean Neg Curv Radii Foliation}
\frac{1}{n} \sum_{i=1}^{n} \frac{2}{1 - \kappa_i} =  1 - \frac
{r^2}{4n(n-1)} S.
\end{equation}
Moreover, if the conformal infinity $(M, [\hg])$ of $(X^{n+1}, g)$
has negative Yamabe invariant, then such foliations are unique.
\end{Cor}

More generally, the relationship (\ref{E: key}) can similarly be applied to the generalized
Yamabe or $\sigma_k$ curvature problem to give foliations of certain hyperbolic Poincar\'{e}
manifolds by hypersurfaces with constant linear combinations or rational  functions of
generalized mean curvatures. For $1 \leq k \leq n$ and
$\lambda = (\lambda_1, \dots, \lambda_n) \in {\R}^n$, let
$$
\sigma_k (\lambda) := \sum_{i_1 < \cdots < i_k} \lambda_{i_1} \cdots
\lambda_{i_k}
$$
denote the $k^{th}$ elementary symmetric function on ${\R}^n$. Let $\Gamma_k$ denote
the connected component of
$$
\{ \lambda \in {\R}^n \ | \ \sigma_k(\lambda) > 0 \}
$$
containing the positive cone $\{ \lambda \in {\R}^n \ | \ \lambda_1,
\dots, \lambda_n > 0 \}$. Given a representative $g_0$ of the
conformal infinity $(M^n, [\hg])$ of a hyperbolic Poincar\'{e}
manifold $(X^{n+1}, g)$, we denote the eigenvalues $\lambda =
(\lambda_1, \dots, \lambda_n)$ of the the Schouten tensor $P_{g_0}$
by $\lambda(P_{g_0})$ and the $k^{th}$ elementary symmetric function
of the eigenvalues of the Schouten tensor $P_{g_0}$ by
$\sigma_k(P_{g_0})$. Moreover, if ${\tilde g}_0 = e^{2\phi_0}g_0$ is
a conformally related metric on $M$, then we denote the $k^{th}$
elementary symmetric function of the eigenvalues of the Schouten
tensor $P_{{\tilde g}_0}$ corresponding to ${\tilde g}_0$ by
$\sigma_k(P_{{\tilde g}_0})$. Applying the works of \cite{GV}
\cite{GW} \cite{LiLi}, it follows from that fact that $M$ is compact
and locally conformally flat, that for $n \geq 3$, if
$\lambda(P_{g_0}) \in \Gamma_k$, then there exists a smooth positive
function $\phi_0$ on $M$ such that ${\tilde g}_0 = e^{2\phi_0}g_0$
with
\begin{equation} \label{E: Conf Sol Gen Yamabe}
\lambda(P_{{\tilde g}_0}) \in \Gamma_k \quad \text{and} \quad \sigma_k(P_{{\tilde g}_0}) = 1.
\end{equation}
In light of (\ref{E: key}), (\ref{E: Conf Sol Gen Yamabe}) and the observations above, we have the following Corollary.

\begin{Cor}\label{C: Gen Cor}
For $n\geq 3$, let $(X^{n+1}, g)$ be a hyperbolic Poincar\'{e} manifold with
conformal infinity $(M^n, [\hg])$. Suppose that there exists a metric $g_0 \in [\hg]$ with
$\lambda(P_{g_0}) \in \Gamma_k$ for some $1 \leq k \leq n$. Then there exists a foliation near $M$
parameterized by a geodesic defining function $r$ associated to a conformal metric
${\tilde g}_0 \in [g_0]$ with constant $\sigma_k$ curvature $\sigma_k(P_{{\tilde g}_0})=1$
such that the level sets of $r$ have outward principal curvatures $\kappa_i=\kappa_i(r)$ satisfying
\begin{equation} \label{E: Gen Foliation}
\sum_{i_1 < \cdots < i_k} \frac{1 + \kappa_{i_1}}{1 -
\kappa_{i_1}}\cdot \frac{1 + \kappa_{i_2}}{1 - \kappa_{i_2}} \cdots
\frac{1 + \kappa_{i_k}}{1 - \kappa_{i_k}} = \l(\frac{r^2}{2}\r)^k.
\end{equation}
\end{Cor}

This paper is organized as follows. In Section \ref{S: AH Mnflds} we introduce hyperbolic Poincar\'{e}
manifolds and we recall several related geometric preliminaries and concepts. In addition, we recall
an application of the ambient metric construction of Fefferman and Graham \cite{FG2}, which gives the
asymptotic expansion (\ref{E: FG}) for the tangential component of hyperbolic Poincar\'{e} metrics in
normal form (\ref{E: FG Norm Form}). In Section \ref{S: Hor Metrics} we introduce the notion of the
horospherical metric associated to a horospherical ovaloid in ${\H}^{n+1}$ and we relate such
horospherical metrics to representatives of the conformal infinity. This observation allows us to put
the two constructions in \cite{JE} and \cite{FG1} \cite{FG2} in the same light. Finally, in Section
\ref{S: Princ Curvs} we prove Theorem \ref{T: Main Thm}.

\section{Hyperbolic Poincar\'{e} Manifolds} \label{S: AH Mnflds}

In this section we introduce hyperbolic Poincar\'{e} manifolds and
their properties mostly adopted from \cite{FG2}. Readers are
referred to \cite{FG2} for details. Let $X^{n+1}$ denote the
interior of a smooth compact manifold ${\bX}^{n+1}$ with boundary
$\partial X = M^n$. A smooth function $r: \bX \to \R$ is said to be
a {\em defining function} for $M$ if

\begin{enumerate}
\item  $r > 0$ in $X$;
\item  $r = 0$ on $M$; and
\item  $d r \neq 0$ on $M$.
\end{enumerate}

A Riemannian metric $g$ on $X$ is then said to be {\em conformally
compact} if for a defining function $r$ for $M$, the conformal
metric $\bg = r^2 g$ extends to a metric on $\bX$. The metric $\bg$
restricted to $T M$ induces a metric $\hg$ on $M$, which rescales
by conformal factor upon change in defining function and therefore
defines a conformal structure $(M, [\hg])$ on $M$ called the
{\em conformal infinity} of $(X, g)$. A straightforward computation
as in \cite{Mazzeo} shows that the sectional curvatures of $g$
approach $-|d r|^2_{\bg}$ near $M$. Accordingly, we have the
following definition for asymptotically hyperbolic manifolds.

\begin{Def} A complete Riemannian manifold  $(X^{n+1}, g)$ is said to be
{\em asymptotically hyperbolic} if $g$ is conformally compact and
$|d r|^2_{\bg}=1$ on $M$ for a defining function $r$ for $M$ in $X$.
\end{Def}

\noindent
We recall the following lemma from \cite{GrahamLee} \cite{Lee} concerning
geodesic defining functions.

\begin{Lem} \label{L: Geod Def Fnct}
Let $(X, g)$ be an asymptotic hyperbolic manifold.  Then any
representative $g_0$ in the conformal infinity of $g$ determines a
unique defining function $r$ such that $r^2 g$ extends to a metric
on $\bX$, $r^2 g |_{TM} = g_0$ and $|dr|^2_{r^2 g} \equiv 1$ in a
neighborhood $U$ of $M$ in $\bX$.
\end{Lem}

We will call such special defining function a geodesic defining
function associated with the representative $g_0$. Given a
representative $g_0$ of the conformal infinity $(M^n, [\hg])$ of an
asymptotic hyperbolic manifold $(X^{n+1}, g)$, the product structure
$M \times [0, \epsilon)$ in a neighborhood of $M$ induced by the
geodesic defining function $r$ from Lemma \ref{L: Geod Def Fnct}
yields the {\em normal form}
$$
g = r^{-2} ( dr^2 + g_r )
$$
with formal asymptotic expansion
$$
g_r = g_0 + r g_1 + r^2 g_2 + \cdots
$$
where the coefficients $g_j$ are symmetric $2$-tensors on $M$.
Decomposing the Einstein tensor $Ric_{g} + ng$ with respect to the
product structure $M \times [0, \epsilon)$ as in \cite{Graham}
yields differential equations that can be successively
differentiated and inductively solved at $r=0$ to derive the
expansions for $n$ odd
\begin{equation} \label{E: G Exp n Odd}
g_r = g_0 + r^2 g_2 + (\text{even powers}) + r^{n-1} g_{n-1} + r^n g_n + \cdots
\end{equation}
while for $n$ even
\begin{equation}  \label{E: G Exp n Even}
g_r = g_0 + r^2 g_2 + (\text{even powers}) + h r^{n} \log r + r^n g_n + \cdots
\end{equation}
provided sufficient regularity is assumed and $Ric_{g} + ng$
vanishes to sufficiently high order at infinity.

For $0 \leq j < n$, the terms $g_j$ in the expansions (\ref{E: G Exp
n Odd}) and (\ref{E: G Exp n Even}) are tensors on $M$ that are
locally determined by the particular representative $g_0$ of the
conformal infinity. For $n$ odd $g_n$ is trace-free but formally
undetermined and for $n$ even $h$ is locally determined and
trace-free while the trace of $g_n$ is locally determined but the
trace-free part of $g_n$ is formally undetermined (see
\cite{Graham}). One can explicitly compute the tensors $g_j$ for $0
\leq j < n$ in the expansions (\ref{E: G Exp n Odd}) and (\ref{E: G
Exp n Even}) using the aforementioned differential equation
resulting from the Einstein condition at infinity. Of particular
interest, for $n\geq3$ one finds that $g_2 = -P_{g_0}$ where
$$
P_{g_0} = \frac{1}{n-2}\l(Ric_{g_0} - \frac{R_{g_0}}{2(n-1)} g_0 \r)
$$
is the Schouten tensor of the conformal representative $g_0$. The
asymptotic expansions described above are fundamental in many works
concerning the geometry and topology of conformally compact
manifolds as well as  in the exploration of properties of
submanifold observables in the AdS/CFT corresponce (see for example
\cite{Anderson} \cite{B-M-Q} \cite{C-Q-Y} \cite{Graham}
\cite{GrahamWitten} \cite{Q}) .

In this note we focus on a class of manifolds that serve as
the prototypical models of asymptotically hyperbolic manifolds known
as {\em hyperbolic Poincar\'{e} manifolds}. Such manifolds are
conformally compact hyperbolic manifolds obtained from quotients of
hyperbolic space ${\H}^{n+1}$ by discrete groups of isometries.
Similar to \cite{Graham}, given a representative $g_0$ of the
conformal infinity $(M^n, [\hg])$ of a hyperbolic Poincar\'{e}
manifold $(X^{n+1}, g)$ one can use the fact that $(X, g)$ has constant
sectional curvature $K_g=-1$ to decompose the tensor
$$
R_{\alpha \beta \gamma \mu} + (g_{\alpha \gamma}g_{\beta \mu} -
g_{\alpha \mu}g_{\beta \gamma}) = 0
$$
with respect to the product structure $M \times [0, \epsilon)$
induced by the geodesic defining function $r$ to yield the
differential equation
\begin{equation} \label{E: Hyp Tang Decomp}
0 = rR^{g^r}_{ijkl}  - \frac{1}{2}(g_{il}g'_{jk} +
g'_{il}g_{jk} - g'_{ik}g_{jl} - g_{ik}g'_{jl}) +
\frac{r}{4}(g'_{il}g'_{jk} - g'_{ik} g'_{jl})
\end{equation}
where latin letters denote tangential directions to $r$ level sets,
$g_{ij} = (g_r)_{ij}$,  and $g'_{ij} = \partial_r g^r_{ij}$ for
simplicity here. Taking successive derivatives of equation (\ref{E:
Hyp Tang Decomp}) and solving at $r=0$ one finds that the tangential
metric
\begin{equation} \label{E: Hyp Exp}
g_r = g_0 - r^2 P_{g_0} + \frac{r^4}{4}Q(P_{g_0})
\end{equation}
where
$$
Q(P_{g_0})_{ij} = g_0^{kl} (P_{g_0})_{ik}(P_{g_0})_{jl}.
$$
The asymptotic expansion (\ref{E: Hyp Exp}) for a hyperbolic
Poincar\'{e} metric is perhaps easier to recognize using the ambient
metric construction of Fefferman and Graham \cite{FG2}. We summarize
the application of the work \cite{FG2} to derive the expansion
(\ref{E: Hyp Exp}) for a hyperbolic Poincar\'{e} metric below.

Let $g_0$ be a representative of the conformal infinity $(M^n,
[\hg])$ of a hyperbolic Poincar\'{e} manifold $(X^{n+1}, g)$ and let
$r$ be the geodesic defining function associated to $g_0$ so that $g$ has the {\em
normal form}
$$
g = r^{-2} ( dr^2 + g_r )
$$
in a neighborhood $M \times [0, \epsilon)$ of $M$. Consider the
ambient metric
\begin{equation} \label{E: Amb Metric}
\tilde{g} = s^2 g - ds^2
\end{equation}
on $M \times [0, \epsilon) \times R_+$. Then $(X, g)$ is
isometrically identified with $\{s = 1\}$ in the ambient spacetime
and a straightforward calculation shows that the curvature tensor of
the ambient metric $\tilde{g}$ satisfies
$$
Riem_{\tilde{g}} = s^2[Riem_g + g \wedge g]
$$
where $(g \wedge g)_{\alpha \beta \gamma \mu} = g_{\alpha
\gamma}g_{\beta \mu} - g_{\alpha \mu}g_{\beta \gamma}$. Hence, it
follows that the ambient metric (\ref{E: Amb Metric}) of a
hyperbolic Poincar\'{e} metric is necessarily flat. Under the change
of variables
$$
-2\rho = r^2, \quad s=rt \quad \text{for} \quad \rho \leq 0
$$
the ambient metric (\ref{E: Amb Metric}) takes the {\em normal form}
$$
\tilde{g} = 2 \rho dt^2 + 2t dt d\rho + t^2 g_{\rho}
$$
where $g_{\rho}$ is a $1$-parameter family of metrics on $M$.
Straightforward computations give the equations
$$
\tilde{R}_{ijkl} = t^2 [R^{g_{\rho}}_{ijkl} +
\frac{1}{2}(g_{il}g'_{jk} + g_{jk}g'_{il} - g_{ik}g'_{jl} -
g_{jl}g'_{ik}) + \frac{\rho}{2}(g'_{ik}g'_{il} - g'_{il}g'_{jk})]
$$
and
$$
\tilde{R}_{\rho ik \rho} = \frac{1}{2} t^2 [g''_{ik} -
\frac{1}{2}g^{jl}g'_{ij}g'_{kl}]
$$
where $g_{ij} = (g_\rho)_{ij}$, $g'_{ij} = \partial_{\rho}
(g_{\rho})_{ij}$ and $g''_{ij} =
\partial_{\rho} \partial_{\rho} (g_{\rho})_{ij}$ for simplicity. Therefore, we
may derive as in the proof of Theorem 7.4 in \cite{FG2} that
$$
g'_{ik}\big|_{\rho = 0} = 2 P^{g_0}_{ik}
$$
and
$$
g''_{ik}\big|_{\rho = 0} = 2 g_0^{jl} P^{g_0}_{ij} P^{g_0}_{kl}
$$
and $g_{ij}'''=0$.

To illustrate the above notions and definitions, we consider the
hyperboloid model of hyperbolic space $(\H^{n+1}, \gH)$. Here
$$
\H^{n+1} = \{ (x, t) \in \R^{n+1,1} \ | \ |x|^2 -t^2 = -1 , \ t > 0
\}
$$
is realized as a hypersurface in Minkowski spacetime $\R^{n+1, 1}$
equipped with the Lorentz metric
$$
g_{_{\L}} = -dt^2 + |dx|^2.
$$
The hyperbolic metric is then given by
$$
\gH = \frac{1}{1+|x|^2} (d|x|)^2 + |x|^2 g_0,
$$
where $g_0$ is the standard round metric on ${\SS}^n$. Letting
$d_{\gH}$ denote the hyperbolic geodesic distance from the vertex
$e_{n+2} \in \H^{n+1} \subset \R^{n+1, 1}$, the function
$$
r := 2e^{-d_{\gH}} = \frac{2}{|x|+\sqrt{1+|x|^2}}
$$
determines the geodesic defining function associated with the
standard round metric as a representative of the conformal infinity
$({\SS}^n, [g_0])$ of hyperbolic space $(\H^{n+1}, \gH)$. We then
have the metric expansion
$$
\gH = r^{-2}\l(dr^2 + \l(1-\frac{r^2}{4}\r)^2 g_0\r).
$$
Notice that $P_{g_0} = \frac{1}{2} g_0$ for the standard round
sphere so that the expansion above is of the form (\ref{E: Hyp
Exp}).

\section{Horospherical Metrics} \label{S: Hor Metrics}

In this section we introduce the horospherical metric on the
space of all horospheres as a parametrization of a neighborhood of
the infinity of hyperbolic space and we present the induced horospherical
metrics on horospherial ovaloids in ${\H}^{n+1}$. Readers are referred
to the paper \cite{JE} for more details. Our goal is to relate
horospherical metrics to representatives of the conformal
infinity and to put the two constructions in \cite{JE} and
\cite{FG1} \cite{FG2} in the same light.

Consider the hyperboloid model of hyperbolic space
$$
\H^{n+1} = \{ (x, t) \in \R^{n+1,1} \ | \ |x|^2 -t^2 = -1 , \ t > 0
\},
$$
where ${\R}^{n+1,1}$ denotes Minkowski spacetime equipped with the
Lorentz metric $g_{_{\L}} = -|dt|^2 + |dx|^2$. Horospheres in
${\H}^{n+1}$ are intersections of degenerate affine hyperplanes of $
{\R}^{n+1,1}$ with ${\H}^{n+1}$ and can be uniquely characterized by
their points at infinity $x \in {\SS}^n$, which are the null
directions inside the hyperplanes, and the signed hyperbolic
distance $\alpha$ from the horosphere to the vertex $e_{n+2} \in
{\H}^{n+1}$, where $\alpha$ is positive if $e_{n+2}$ is inside a
given horosphere and negative otherwise. Accordingly, one can
identify the space of horospheres in ${\H}^{n+1}$ with ${\SS}^n
\times {\R}$ and endow the space of horospheres with a natural
degenerate metric $\langle \cdot \ ,\cdot \rangle_{\infty} =
e^{2\alpha} g_0$ in the conformal class of the round metric $g_0$ on
${\SS}^n$.

Now suppose
$$
\phi: \Sigma^n \to {\H}^{n+1}
$$
is an immersed oriented hypersurface and let
$$
\eta: {\Sigma}^n \to {\SS}^{n+1}_1
$$
denote the Lorentzian unit normal map taking values in de-Sitter
spacetime
$$
{\SS}^{n+1}_1 = \{ (x, t) \in \R^{n+1,1} \ | \ |x|^2 -t^2 = 1\}.
$$
From the map
\begin{equation}\label{E: HO Light Cone Map}
\psi:= \phi + \eta: {\Sigma}^n \to {\N}^{n+1}_+
\end{equation}
taking values in the positive light-cone
$$
{\N}^{n+1}_+ = \{ (x, t) \in \R^{n+1,1} \ | \ |x|^2 -t^2 = 0, \ t >
0\},
$$
one defines the hyperbolic Gauss map as the direction of the
light-cone map (\ref{E: HO Light Cone Map}) in ${\SS}^n$.  One finds
that the light-cone map of horospheres is constant for the inward
orientation and that parallel horospheres correspond to collinear
vectors in ${\N}^{n+1}_+$. Hence, one also can identify the space of
horospheres in ${\H}^{n+1}$ with ${\N}^{n+1}_+$. Moreover, it is
easily seen that the horospherical metric on the space of all
horospheres is exactly the same as the induced metric on the
light-cone from the Lorentz metric $g_{_{\L}}$ of Minkowski
spacetime.

One therefore can realize the {\em horospherical metric} associated
to a horospherical ovaloid in ${\H}^{n+1}$, that is a compact
hypersurface $\Sigma^n \subseteq {\H}^{n+1}$ for which the Gauss map
is regular, as the pull-back by the light-cone map $\psi$ of the
induced metric on the hypersurface as viewed in the positive
light-cone ${\N}^{n+1}_+$. We recall from \cite{JE} that a compact
immersed hypersurface is said to be a horospherical ovaloid in
$\H^{n+1}$ if it can be oriented so that it is horospherically
convex at every point and that an oriented hypersurface in
$\H^{n+1}$ is horospherically convex at a point if and only if all
the principal curvatures of at the point verify simultaneously less
than $1$ or greater than $1$.

Alternatively, one can define the horospherical metric as in
\cite{Schlenker} by
\begin{equation} \label{E: Horo Metric}
h:= I_{\Sigma} -  2 II_{\Sigma} + III_{\Sigma}
\end{equation}
where $I_{\Sigma}, II_{\Sigma}$ and $III_{\Sigma}$ are respectively
the first, second and third fundamental forms of $\Sigma$ in
${\H}^{n+1}$. In \cite{JE} Espinar, G\'{a}lvez and Mira view the
image of the light cone map (\ref{E: HO Light Cone Map}) as a
co-dimension $2$ hypersurface in Minkowoski spacetime and derive a
relation between the principal curvatures of an immersed
hypersurface in ${\H}^{n+1}$ and the eigenvalues of the Schouten
tensor of its associated horospherical metric. In order to connect
the work of \cite{JE} with ours in the context of conformal
geometry, we compute the horospherical metric as defined in (\ref{E:
Horo Metric}). Given a hyperbolic Poincar\'{e} manifold $(X^{n+1},
g)$ and a respresentative $\gamma$ of its conformal infinity $(M^n,
[\hg])$ we first compute the third fundamental form on level sets
determined by the associated geodesic defining function $r$. A
straightforward computation gives
$$
\aligned &III_r(\partial_{i}, \partial_{j}) =
I_r(\nabla_{\partial_i} N_r, \nabla_{\partial_j} N_r) =
I_r(\nabla_{\partial_i} r
\partial_r,
\nabla_{\partial_j} r \partial_r) \notag \\
&= {r}^{-2}g^r_{ij} - {r}^{-1} \partial_r g^r_{ij} + \frac{1}{4}
g_r^{pq} \partial_r g^r_{ip} \partial_r g^r_{jq}\notag\\
&=\l(r^{-2} \gamma_{ij} - P^{\gamma}_{ij} + \frac{r^2}{4}
\gamma^{kl} P^{\gamma}_{ik} P^{\gamma}_{jl}\r) + \l(2P^{\gamma}_{ij} -
r^2 \gamma^{kl}  P^{\gamma}_{ik} P^{\gamma}_{jl}\r)\notag\\
&+ \frac{1}{4}g_r^{pq}\l(-2r P^{\gamma}_{ip} + r^3 \gamma^{kl}
P^{\gamma}_{ik} P^{\gamma}_{pl}\r)\l(-2r P^{\gamma}_{jq} + r^3 \gamma^{kl}
P^{\gamma}_{jk} P^{\gamma}_{ql}\r).
\endaligned
$$
In terms of an orthonormal basis $\{e_1, \dots, e_n\}$ with respect
to $\gamma$ that diagonalizes the tensor $P_{\gamma}$, it follows
$$
\aligned III_r(e_i, e_j) &=\l(r^{-2} \delta_{ij} - \lambda_{i}
\delta_{ij} + \frac{r^2}{4} \lambda_{i}^{2} \delta_{ij}\r) +
\l(2\lambda_{i} \delta_{ij} -
r^2 \lambda_{i}^{2} \delta_{ij}\r)\notag \\
&+ \frac{1}{4} \l(1 - \frac{r^2}{2} \lambda_{k}\r)^{-2} \delta^{kl}
\l(-2r \lambda_{i} \delta_{ik} + r^3 \lambda_{i}^2 \delta_{ik}\r)
\l(-2r \lambda_{j}\delta_{jl} + r^3 \lambda_{j}^{2}\delta_{jl}\r)\notag \\
&=r^{-2}  \l(1 + \frac{r^2}{2} \lambda_{i}\r)^2  \delta_{ij}.
\endaligned
$$
Therefore, the horospherical metric associated to a level set of a geodesic
defining function $r$ is
$$
\aligned
h(e_i, e_j) & = I_r(e_i, e_j)  - 2II_r(e_i, e_j) +III_r(e_i, e_j)\notag \\
&=r^{-2} \l(1 - \frac{r^2}{2} \lambda_{i}\r)^2\delta_{ij} +2r^{-2}
\l(1 - \frac{r^2}{2} \lambda_{i}\r)  \l(1 + \frac{r^2}{2}
\lambda_{i}\r) \delta_{ij} \notag\\ & \quad +r^{-2}  \l(1 +
\frac{r^2}{2} \lambda_{i}\r)^2
\delta_{ij}\notag \\
&=r^{-2} \l(1 - r^2 \lambda_{i} +  \frac{r^4}{4} \lambda_{i}^{2} +2 -
\frac{r^4}{2} \lambda_{i}^{2} + 1 + r^2 \lambda_{i} + \frac{r^4}{4}
\lambda_{i}^2\r) \delta_{ij}\notag \\
&=4 r^{-2} \delta_{ij}.
\endaligned
$$
Thus, given a conformal representative $\gamma$ of the conformal
infinity $(M^n, [\hg])$ of a hyperbolic Poincar\'{e} manifold
$(X^{n+1}, g)$ and its associated geodesic defining function $r$,
the horospherical metrics associated to the level sets of $r$ are given by $h=4
r^{-2}\gamma$. On the other hand, given an outwardly convex smooth
hypersurface $\Sigma^n \subset X$, from which the exponential map
is a diffeomorphism from the normal bundle to the outside, we find from the
associated geodesic defining function $\tilde{r} = e^{-d_{\Sigma}}$,
that the horospherical metric on $\Sigma = \{\tilde{r} =1\}$
is given by $h \in [\hg]$. Hence, in the context of conformal
geometry one may regard horospherical metrics associated to the
hypersurfaces given in \cite{JE} simply as conformal representatives of the
conformal infinity.

Now we would like to illustrate that the ambient metric construction
in \cite{FG2} somehow gives a nice extension to the notions of the horospherical
metrics in \cite{JE}. As in \cite{FG2}, given a hyperbolic
Poincar\'{e} manifold $(X^{n+1}, g)$ with conformal infinity $(M,
[\hg])$ we define the metric bundle ${\mathcal G}$ over $M$ to be
the space of pairs $(h, x)$ with $x \in M$ and $h=s^2 \hg(x)$ for
some $s \in {\R}_+$ where ${\mathcal G}$ is equipped with the
projection
$$
\pi : {\mathcal G} \to M \quad \text{defined by} \quad (h, x)
\overset{\pi}\mapsto x
$$
and dilations
$$
\delta_s : {\mathcal G} \to {\mathcal G} \quad \text{defined by}
\quad (h, x) \overset{\delta_s}\mapsto (s^2 h, x)
$$
for $s \in {\R}_+$. The metric bundle ${\mathcal G}$ assumes the
role of the light cone, that is, the space of all horospheres, and
the metric bundle is similarly equipped with a tautological
degenerate metric defined at $z=(h, x) \in {\mathcal G}$ by
$g_{_{\mathcal G}} = \pi^* h$, which is homogeneous of degree $2$ with
respect to dilations and therefore depends only on the conformal
class $[\hg]$.

Fixing a representative $g_0$ of the conformal infinity $(M,
[\hg])$, one obtains a trivialization of metric bundle ${\mathcal G}
\cong {\R}_+ \times M$ by identifying
$$
(t, x) \in {\R}_+ \times M \quad \text{with} \quad (t^2 g_0(x), x)
\in {\mathcal G}.
$$
Given local coordinates $(x) = (x^1, \dots, x^n)$ on ${\mathcal U}
\subset M$ we obtain local coordinates $(t, x)$ on
$\pi^{-1}({\mathcal U})$ where
$$
g^{\mathcal G}_{ij} = t^2 g^{0}_{ij} dx^i dx^j
$$
so that the representative $g_0$ of the conformal infinity of $(X,
g)$ can be considered as the section of the bundle ${\mathcal G}
\cong {\R}_+ \times M$ determined by the level submanifold
$\{t=1\}$. On the ambient space $({\R}_+ \times M) \times {\R}$ with
coordinates $(t, x, \rho)$ the ambient or cone metric $\tilde{g} =
s^2 g - ds^2$ from (\ref{E: Amb Metric}) with $(X, g) = \{s = 1\}$
takes the normal form
$$
\tilde{g} = 2 \rho dt^2 + 2t dt d\rho + t^2 g_{\rho}
$$
where
$$
-2\rho = r^2, \quad s=rt \quad \text{for} \quad \rho \leq 0
$$
and $r$ is the geodesic defining function uniquely associated
to $g_0$. Therefore, given an outwardly convex hypersurface
$\Sigma^n \subset X^{n+1}$ and letting $\alpha = d_{\Sigma}$
denote the signed geodesic distance from $\Sigma$, which is
positive outside $\Sigma$, one finds that under the change of
variables
$$
t=e^{\alpha}
$$
that the ambient metric restricted to $(X, g) = \{s = 1\}$
takes the form
$$
\tilde{g} \big|_{X} = d\alpha^2 + e^{2\alpha} g_{\alpha}.
$$
Hence, one may view the change of variables $t=e^{\alpha}$ with
respect to a given hypersurface as straightening out the hypersurface
and giving a new coordinate on the metric bundle, which results in
determining a new representative of the conformal infinity.

\section{Principal Curvatures} \label{S: Princ Curvs}

In this section we carry out a straightforward calculation to prove our main
theorem. Suppose that $(X^{n+1}, g)$ is a hyperbolic Poincar\'{e}
manifold and $(M^n, [\hg])$ is its conformal infinity. Let
$\gamma$ be a representative of the conformal infinity and let $r$
be the geodesic defining function associated to $\gamma$ so that
$g$ has the normal form
$$
g = r^{-2} (dr^2 + g_r)
$$
near $M$ with
$$
g_r = \gamma -  r^2 P_{\gamma} + \frac{r^4}{4} Q(P_\gamma)
$$
where
$$
Q(P_\gamma)_{ij} = \gamma^{kl}(P_\gamma)_{ik}(P_\gamma)_{jl}.
$$
Then the level sets of $r$ give a foliation near $M$ with induced
metric
$$
I_r = r^{-2} g_r = r^{-2} \gamma - P_{\gamma} +
 \frac{r^2}{4} Q(P_{\gamma})
$$
and outward pointing normal $N_r = -r \partial_r$ where $\partial_r
:= \nabla_{\bg} r$. Hence, the level sets of $r$ have second
fundamental form, according to our definition
(\ref{secondfundamentalform}),
$$
\aligned II_r &= \frac{1}{2} r \partial_r \l( r^{-2} g_r \r) =
- r^{-2} g_r + \frac{1}{2} r^{-1} \partial_r g_r\notag \\
& = - r^{-2} \gamma + P_{\gamma} - \frac{r^2}{4}
Q(P_{\gamma}) - P_{\gamma} + \frac{r^2}{2} Q(P_{\gamma})\notag \\
&=  - r^{-2} \gamma + \frac{r^2}{4} Q(P_{\gamma}).
\endaligned
$$
Now let $\{e_1, \dots, e_n\}$ denote an orthonormal basis with
respect to $\gamma$ that diagonalizes the tensor $P_{\gamma}$. Then
$$
\gamma(e_i, e_j) = \delta_{ij} \quad \text{and} \quad
P_{\gamma}(e_i, e_j) = \lambda_i \delta_{ij}
$$
where $\lambda_i$ denotes the $i^{th}$ eigenvalue of the tensor
$P_{\gamma}$. Moreover,
$$
\aligned I_r (e_i, e_j) =r^{-2} \gamma(e_i, e_j)  &- P_{\gamma}(e_i,
e_j) + \frac{r^2}{4}  \gamma^{-1}(e_k, e_l) P_{\gamma}(e_i, e_k)
P_{\gamma}(e_j, e_l)\notag \\
&=r^{-2} \delta_{ij} - \lambda_i \delta_{ij} + \frac{r^2}{4}
\delta^{kl} \lambda_i \delta_{ik}  \lambda_j \delta_{jl}\notag \\
&=r^{-2} \l(1 - r^{2} \lambda_i + \frac{r^4}{4} \lambda_i^2\r)
\delta_{ij}\notag \\
&=r^{-2} \l(1 - \frac{r^2}{2} \lambda_i\r)^2 \delta_{ij}
\endaligned
$$
and
$$
\aligned II_r (e_i, e_j) = - r^{-2} &\gamma(e_i, e_j)  +
\frac{r^2}{4} \gamma^{-1}(e_k, e_l) P_{\gamma}(e_i, e_k)
P_{\gamma}(e_j, e_l)
\notag \\
&= - r^{-2} \delta_{ij} +  \frac{r^2}{4}  \delta^{kl} \lambda_i
\delta_{ik}  \lambda_j \delta_{jl}\notag \\
&= - r^{-2} \l(1 - \frac{r^4}{4} \lambda_i^2\r) \delta_{ij}
\notag \\
&= - r^{-2} \l(1 - \frac{r^2}{2} \lambda_i\r) \l(1 + \frac{r^2}{2}
\lambda_i\r) \delta_{ij}.
\endaligned
$$
Therefore,
$$
\aligned &\l(I_r^{-1}II_r\r) (e_i, e_j) = I_r^{-1} (e_i, e_k) II_r
(e_k, e_j)\notag \\
&= -r^{2} \l(1 - \frac{r^2}{2} \lambda_i\r)^{-2} \delta^{ik} r^{-2}
\l(1 - \frac{r^2}{2} \lambda_k\r) \l(1 + \frac{r^2}{2} \lambda_k\r)
\delta_{kj} =  - \frac{1 + \frac{r^2}{2} \lambda_i} {1 -
\frac{r^2}{2} \lambda_i} \delta_{ij}.
\endaligned
$$
But the Weingarten matrix $I_r^{-1}II_r$ on the level sets of $r$
satisfies
$$
\l(I_r^{-1}II_r\r) (e_i, e_j) = \kappa_i^r \delta_{ij}
$$
where $\kappa_i^r$ denotes the $i^{th}$ principal curvature of a
level set of $r$ with respect to the outward direction. Hence,
$$
\kappa_i^r = - \frac{1 + \frac{r^2}{2} \lambda_i}{1 - \frac{r^2}{2}
\lambda_i} =  - \frac{2}{1 - \frac{r^2}{2} \lambda_i} + 1
$$
so that
$$
1 - \frac{r^2}{2} \lambda_i = \frac{2}{1 - \kappa_i^r},
$$
which establishes Theorem \ref{T: Main Thm}.


\begin{thebibliography}{99}

\bibitem{Anderson}
M. Anderson.
\newblock $L^2$ Curvature and Volume Renormalization of AHE Metrics on 4-manifolds.
\newblock {\em Math. Res. Lett.} 8: 171-188, 2001.


\bibitem{BM}
Eric Bahuaud and Tracey Marsh.
\newblock Holder Compactification for Some Manifolds with Pinched Negative Curvature at Infinity.
\newblock arXiv:0601503v1 to appear {\em Canadian Journal of Mathematics}.

\bibitem{B-M-Q}
Vincent Bonini, Pengzi Miao and Jie Qing.
\newblock Ricci Curvature Rigidity for Weakly Asymptotically Hyperbolic Manifolds.
\newblock {\em Comm. Anal. Geom.} 14(3): 603-612, 2006.

\bibitem{C-Q-Y}
Alice Chang, Jie Qing and Paul Yang.
\newblock On the Topology of Conformally Compact Einstein 4-manifolds.
\newblock {\em math.DG: arXiv:0305085}, Noncompact problems at the intersection of geometry, analysis,
and topology, 49--61, {\em Contemp. Math.}, 350, Amer. Math. Soc.,
Providence, RI, 2004.

\bibitem{Epstein1}
Charles L. Epstein.
\newblock Envelopes of Horospheres and Weingarten Surfaces in Hyperbolic
$3$-Space, unpublished preprint.

\bibitem{Epstein2} Charles L. Epstein.
\newblock An Asymptotic Volume Formula for Convex Cocompact
Hyperbolic Manifolds.
\newblock Appendix A: S. J. Patterson and P. A. Perry.
\newblock The Divisor of Selberg Zeta Function for Kleinian Groups.
\newblock Duke Math. J. 106: 321-390, 2001.

\bibitem{JE}
Jos\'{e} M. Espinar, Jos\'{e} A. G\'{a}lvez and Pablo Mira.
\newblock Hypersurfaces in ${\H}^{n+1}$ and Conformally Invariant Equations: The Generalized Christoffel and Nirenberg Problems.
\newblock {\em Jour. European Math. Soc.}, 11 (2009) 903-939.

\bibitem{Graham}
C.~Robin Graham.
\newblock Volume and Area Renormalizations for Conformally Compact Einstein Metrics.
\newblock {\em math.DG: arXiv:9909042}, The Proceedings of the 19th Winter School
``Geometry and Physics'' (Srn¨ª, 1999). {\em Rend. Circ. Mat.
Palermo} (2) Suppl. No. 63 (2000), 31--42.

\bibitem{FG1}
Charles Fefferman and C.~Robin Graham.
\newblock Conformal Invariants.
\newblock {\em Asterisque} 95-116, 1985.

\bibitem{FG2}
Charles Fefferman and C.~Robin Graham.
\newblock The Ambient Metric.
\newblock {\em math.DG: arXiv:0710.0919}

\bibitem{GrahamLee}
C.~Robin Graham and John~M. Lee.
\newblock Einstein Metrics with Prescribed Conformal Infinity on the Ball.
\newblock {\em Adv. Math.} 87(2): 186--225, 1991.

\bibitem{GrahamWitten}
C.~Robin Graham and Edward Witten.
\newblock Confomal Anomaly of Submanifold Observables in the Ads/CFT Correspondence.
\newblock {\em Nucl.Phys.} B546: 52-64, 1999.

\bibitem{GW}
P. Guan and G. Wang.
\newblock A Fully Nonlinear Conformal Flow on Locally Conformally Flat Manifolds.
\newblock {\em J. Reine Angew. Math.} 557: 219--238, 2003.

\bibitem{GV}
M. Gursky and J. Viaclovsky.
\newblock Prescribing symmetric functions of the eigenvalues of the Ricci
tensor.
\newblock {\em Annals of Mathematics} 166 (2007) no. 2, 475-531.

\bibitem{Lee}
John~M. Lee.
\newblock The Spectrum of an Asymptotically Hyperbolic {E}instein Manifold.
\newblock {\em Comm. Anal. Geom.} 3(1-2): 253--271, 1995.

\bibitem{LiLi}
A. Li and Y.Y. Li.
\newblock On Some Conformally Invariant Fully Nonlinear Equations.
\newblock {\em Comm. Pure Appl. Math.} 56: 1416--1464, 2003.

\bibitem{Mazzeo}
Rafe Mazzeo.
\newblock The {H}odge Cohomology of a Conformally Compact Metric.
\newblock {\em J. Differential Geom.} 28(2): 309--339, 1988.

\bibitem{MP}
Rafe Mazzeo and Frank Pacard.
\newblock Constant Curvature Foliations on Asymptotically Hyperbolic Spaces.
\newblock math.DG: arXiv:0710.2298.

\bibitem{Q}  Jie Qing. \newblock On the Rigidity for Conformally Compact Einstein Manifolds
\newblock {\em Int. Math. Res. Not.} 21: 1141-1153, 2003.

\bibitem{Schlenker}
J. M. Schlenker.
\newblock Hypersurfaces in ${\H}^n$ and the Space of Its Horospheres
\newblock {\em Geom. Anal. Funct. Anal.} 12: 395--435, 2002.
\end{thebibliography}
\end{document}